\documentclass[11pt]{article}

\usepackage{amssymb,latexsym,amsmath}

\usepackage{graphicx}

\begin{document}

\newcommand{\bfi}{\bfseries\itshape}

\makeatletter

\@addtoreset{figure}{section}

\def\thefigure{\thesection.\@arabic\c@figure}

\def\fps@figure{h, t}

\@addtoreset{table}{bsection}

\def\thetable{\thesection.\@arabic\c@table}

\def\fps@table{h, t}

\@addtoreset{equation}{section}

\def\theequation{\thesubsection.\arabic{equation}}

\makeatother

\newtheorem{thm}{Theorem}[section]

\newtheorem{prop}[thm]{Proposition}

\newtheorem{lema}[thm]{Lemma}

\newtheorem{cor}[thm]{Corollary}

\newtheorem{defi}[thm]{Definition}

\newtheorem{rk}[thm]{Remark}

\newtheorem{exempl}{Example}[section]

\newenvironment{exemplu}{\begin{exempl}  \em}{\hfill $\surd$

\end{exempl}}

\newcommand{\comment}[1]{\par\noindent{\raggedright\texttt{#1}

\par\marginpar{\textsc{Comment}}}}

\newcommand{\todo}[1]{\vspace{5 mm}\par \noindent \marginpar{\textsc{ToDo}}\framebox{\begin{minipage}[c]{0.95 \textwidth}

\tt #1 \end{minipage}}\vspace{5 mm}\par}

\newcommand{\ea}{\mbox{{\bf a}}}

\newcommand{\eu}{\mbox{{\bf u}}}

\newcommand{\ueu}{\underline{\eu}}

\newcommand{\ueo}{\overline{u}}

\newcommand{\oeu}{\overline{\eu}}

\newcommand{\ew}{\mbox{{\bf w}}}

\newcommand{\ef}{\mbox{{\bf f}}}

\newcommand{\eF}{\mbox{{\bf F}}}

\newcommand{\eC}{\mbox{{\bf C}}}

\newcommand{\en}{\mbox{{\bf n}}}

\newcommand{\eT}{\mbox{{\bf T}}}

\newcommand{\eL}{\mbox{{\bf L}}}

\newcommand{\eR}{\mbox{{\bf R}}}

\newcommand{\eV}{\mbox{{\bf V}}}

\newcommand{\eU}{\mbox{{\bf U}}}

\newcommand{\ev}{\mbox{{\bf v}}}

\newcommand{\eve}{\mbox{{\bf e}}}

\newcommand{\uev}{\underline{\ev}}

\newcommand{\eY}{\mbox{{\bf Y}}}

\newcommand{\eK}{\mbox{{\bf K}}}

\newcommand{\eP}{\mbox{{\bf P}}}

\newcommand{\eS}{\mbox{{\bf S}}}

\newcommand{\eJ}{\mbox{{\bf J}}}

\newcommand{\eB}{\mbox{{\bf B}}}

\newcommand{\eH}{\mbox{{\bf H}}}

\newcommand{\leb}{\mathcal{ L}^{n}}

\newcommand{\eI}{\mathcal{ I}}

\newcommand{\eE}{\mathcal{ E}}

\newcommand{\hen}{\mathcal{H}^{n-1}}

\newcommand{\eBV}{\mbox{{\bf BV}}}

\newcommand{\eA}{\mbox{{\bf A}}}

\newcommand{\eSBV}{\mbox{{\bf SBV}}}

\newcommand{\eBD}{\mbox{{\bf BD}}}

\newcommand{\eSBD}{\mbox{{\bf SBD}}}

\newcommand{\ecs}{\mbox{{\bf X}}}

\newcommand{\eg}{\mbox{{\bf g}}}

\newcommand{\paromega}{\partial \Omega}

\newcommand{\gau}{\Gamma_{u}}

\newcommand{\gaf}{\Gamma_{f}}

\newcommand{\sig}{{\bf \sigma}}

\newcommand{\gac}{\Gamma_{\mbox{{\bf c}}}}

\newcommand{\deu}{\dot{\eu}}

\newcommand{\dueu}{\underline{\deu}}

\newcommand{\dev}{\dot{\ev}}

\newcommand{\duev}{\underline{\dev}}

\newcommand{\weak}{\stackrel{w}{\approx}}

\newcommand{\mild}{\stackrel{m}{\approx}}

\newcommand{\strong}{\stackrel{s}{\approx}}

\newcommand{\weakdown}{\rightharpoondown}

\newcommand{\opg}{\stackrel{\mathfrak{g}}{\cdot}}

\newcommand{\opunu}{\stackrel{1}{\cdot}}
\newcommand{\opdoi}{\stackrel{2}{\cdot}}

\newcommand{\opn}{\stackrel{\mathfrak{n}}{\cdot}}

\newcommand{\tr}{\ \mbox{tr}}

\newcommand{\Ad}{\ \mbox{Ad}}

\newcommand{\ad}{\ \mbox{ad}}

\renewcommand{\contentsname}{ }

\title{Curvature of sub-Riemannian spaces}

\author{Marius Buliga \\
 \\
Institut Bernoulli\\
B\^{a}timent MA \\
\'Ecole Polytechnique F\'ed\'erale de Lausanne\\
CH 1015 Lausanne, Switzerland\\
{\footnotesize Marius.Buliga@epfl.ch} \\
 \\
\and
and \\
 \\
Institute of Mathematics, Romanian Academy \\
P.O. BOX 1-764, RO 70700\\
Bucure\c sti, Romania\\
{\footnotesize Marius.Buliga@imar.ro}}

\date{This version: 26.11.2003}

\maketitle

\newpage

\section{Introduction}

To any metric spaces there is an associated metric profile. The rectifiability of the metric 
profile gives a good notion of curvature of a sub-Riemannian space. We shall say that a 
curvature class is the rectifiability class of the metric profile. We classify then the curvatures 
by looking to homogeneous metric spaces.  The classification 
problem is solved for contact 3 manifolds, where we rediscover a 3 dimensional family of 
homogeneous contact manifolds, with a distinguished  2 dimensional family of  contact 
manifolds which don't have a natural group structure.  The classification of 3 dimensional 
homogeneous contact manifolds has been done by Hughen \cite{hughen}. 

I have discovered metric profiles in various proof of Mitchell theorem 1. Also this is explained  
in the paper. In my opinion, the use of the notion of metric profile clarifies the question: why several 
proofs for same result (Mitchell theorem 1) and moreover, any of them equally long and complex? 

It has to be mentioned that contrary to other attempts to define the curvature of a sub-Riemannian 
manifold, here is presented an almost pure metrical construction, not using differential geometry, which 
is notoriously misleading when used in a sub-Riemannian frame. Once one knows what to look 
for, then differential geometry (read "Euclidean analytic differential geometry") recovers its well 
known strength, though. 

The structure of the paper is described further. In sections 2 -- 5 is given a short presentation of 
sub-Riemannian manifolds, Carnot groups, Pansu derivative and Gromov-Hausdorff distance. 
For the expert reader these sections serve only to fix notations needed later. 

Section 5 is about deformations of sub-Riemannian manifold, seen as curves in the space CMS of 
isometry classses of compact metric spaces, with the Gromov-Hausdorff distance.  

In section 6 can be found a discussion of various proofs of Mitchell \cite{mit} theorem 1. This section 
justifies the notion of metric profile, which is the subject of section 7. In the same section is given thee notion of curvature in terms of rectifiability classes of metric profiles. 

In order to classify the curvatures homogeneous spaces are used. Section 8 is dedicated to this 
subject. 

As an application, in section 9 are studied the homogeneous contact 3 manifolds, Finally, in section 10 the problem of classification is solved for a large class of contact 3 manifolds.

\section{Regular sub-Riemannian manifolds} 

Classical references to this subject are Bella\"{\i}che \cite{bell} and Gromov \cite{gromo}. The interested reader is advised to look also to the references of these papers.  

Let $M$ be a connected manifold. A distribution (or horizontal bundle) is a subbundle 
$D$ of $M$. To any point $x \in M$ there is associated the vectorspace $D_{x} \subset 
T_{x}M$. 

Given the distribution $D$, a point $x \in M$ and a sufficiently small  open neighbourhood $x \in U \subset M$, one can define on $U$ a filtration of bundles
as follows. Define first the class of horizontal vectorfields on $U$:
$$\mathcal{X}^{1}(U,D) \ = \ \left\{ X \in \Gamma^{\infty}(TU) \mbox{ : }
\forall y \in U \ , \ X(u) \in D_{y} \right\}$$
Next, define inductively for all positive integers $k$:
$$ \mathcal{X}^{k+1} (U,D) \ = \ \mathcal{X}^{k}(U,D) \cup [ \mathcal{X}^{1}(U,D),
\mathcal{X}^{k}(U,D)]$$
Here $[ \cdot , \cdot ]$ denotes vectorfields bracket. We obtain therefore a filtration $\displaystyle \mathcal{X}^{k}(U,D) \subset \mathcal{X}^{k+1} (U,D)$.
Evaluate now this filtration at $x$:
$$V^{k}(x,U,D) \ = \ \left\{ X(x) \mbox{ : } X \in \mathcal{X}^{k}(U,D)\right\}$$
There are
$m(x)$, positive integer, and small enough $U$ such that $\displaystyle
V^{k}(x,U,D) = V^{k}(x,D)$ for all $k \geq m$ and
$$D_{x} V^{1}(x,D) \subset V^{2}(x,D) \subset ... \subset V^{m(x)}(x,D)$$
We equally have
$$ \nu_{1}(x) = \dim V^{1}(x,D) < \nu_{2}(x) = \dim V^{2}(x,D) < ... < n = \dim M$$
Generally $m(x)$, $\nu_{k}(x)$
may vary from a point to another.

The number $m(x)$ is called the step of the distribution at 
$x$.

\begin{defi}
The distribution $D$ is regular if $m(x)$, $\nu_{k}(x)$ are constant on the manifold $M$. 

The distribution is completely non-integrable if for any $x \in M$ we have $\displaystyle 
V^{m(x)} = T_{x}M$. 
\label{dreg}
\end{defi}

\begin{defi}
A sub-Riemannian (SR) manifold is a triple $(M,H, g)$, where $M$ is a
connected manifold, $H$ is a completely non-integrable distribution on $M$, and $g$ is a metric (Euclidean inner-product) on the
horizontal bundle $H$.
\label{defsr}
\end{defi}

The Carnot-Carath\'eodory distance associated to the sub-Riemannian manifold is the 
distance induced by the length $l$ of horizontal curves:
$$d(x,y) \ = \ \inf \left\{ l(c) \mbox{ : } c: [a,b] \rightarrow M \
, \ c(a) = x \ , \  c(b) = y \right\} $$

The Chow theorem ensures the existence of a horizontal path linking any two sufficiently 
closed points, therefore the CC distance it at least locally finite. 

We shall work further only with regular sub-Riemannian manifolds, if not otherwise stated.

Bella\"{\i}che introduced the concept of privileged chart around a point $p \in M$.  

Let $(x_{1}, ... , x_{n}) \mapsto \phi(x_{1}, ... , x_{n}) \in M$ be a chart of $M$ around $p$ (i.e. $p$ has coordinates $(0,....,0)$).  Denote by 
$X_{1}, ... , X_{n}$ the frame of vectorfields  associated to the coordinate chart.  The chart is called adapted  (or the frame is called adapted) if the following happens: $X_{1}, ... , X_{\nu_{1}}$ forms a
basis of $V^{1}$, $X_{\nu_{1}+1}, ... , X_{\nu_{2}}$ form a basis of $V^{2}$, and so on.  

Suppose that the frame $X_{1}, ... , X_{n}$ is adapted. The degree of $X_{i}$ is then $k$ if 
$X_{i} \in  V^{k} \setminus V^{k-1}$.

\begin{defi}
A chart (or a frame) is privileged if moreover the following happens: for any $i = 1, ... , n$ the 
function 
$$ t \mapsto d(p, \phi( ... , t, ...)) $$ 
(with $t$ on the position $i$) is exactly of order  $deg \  X_{i}$  at $t=0$. 
\end{defi}

Privileged charts (frames) always exist, as proved by Bella\"{\i}che  \cite{bell} Theorem 4.15. 

A privileged frame transforms the filtration into a direct sum. Define 
$$V_{i} = \ span \ \left\{ X_{k} \mbox{ : } deg \ X_{k}  = i  \right\}$$
Then the tangent  space decomposes as a direct sum of vectorspaces $V_{i}$. Moreover, each 
space $V^{i}$ decomposes in a direct sum of spaces $V_{k}$ with $k \leq i$. 

The intrinsic dilatations associated to a privileged frame are defined, in the chart 
$\phi$,  for any $\varepsilon > 0$ (sufficiently small if necessary)  by 
$$\delta_{\varepsilon} (x_{i})  = (\varepsilon^{deg \ i} x_{i} )$$
We may define (locally around $p$) a Lie bracket associated to the privileged frame, which comes from the vectorfield bracket written in coordinates with respect to the  frame (which is a basis of the tangent space). 

In terms of vectorfields, the intrinsic dilatation associated to the privileged frame transforms 
$X_{i}$ into $$\Delta_{\varepsilon} X_{i} = \varepsilon^{deg \ X_{i}} X_{i}$$ and the metric $g$ into $\frac{1}{\varepsilon^{2}} g$. 

The nilpotentization of the distribution with respect to the chosen privileged frame is then 
the bracket 
$$[X,Y]_{N} = \lim_{\varepsilon \rightarrow 0} \Delta_{\varepsilon}^{-1} [\Delta_{\varepsilon} X , 
\Delta_{\varepsilon}Y]$$ 

It is very important to notice that the useful part of the nilpotentization bracket  is its evaluation at the point $p$. It is generically false that there are privileged coordinates around an open set in $M$. This is however true in the particular case of contact manifolds, as a consequence of Frobenius theorem.

\section{Carnot groups}

Carnot groups are particular examples of sub-Riemannian manifolds. They are especially important because they provide infinitesimal models for any sub-Riemannian manifold.

\begin{defi}
A Carnot (or stratified nilpotent) group is a
connected simply connected group $N$  with  a distinguished vectorspace
$V_{1}$ such that the Lie algebra of the group has the
direct sum decomposition:
$$n \ = \ \sum_{i=1}^{m} V_{i} \ , \ \ V_{i+1} \ = \ [V_{1},V_{i}]$$
The number $m$ is the step of the group. The number
$$Q \ = \ \sum_{i=1}^{m} i \ dim V_{i}$$
is called the homogeneous dimension of the group.
\label{dccgroup}
\end{defi}

Because the group is nilpotent and simply connected, the
exponential mapping is a diffeomorphism. We shall identify the
group with the algebra, if is not locally otherwise stated.

The structure that we obtain is a set $N$ endowed with a Lie
bracket and a group multiplication operation given by the Baker-Campbell-Hausdorff formula. 

Any Carnot group admits a one-parameter family of dilatations. For any
$\varepsilon > 0$, the associated dilatation is:
$$ x \ = \ \sum_{i=1}^{m} x_{i} \ \mapsto \ \delta_{\varepsilon} x \
= \ \sum_{i=1}^{m} \varepsilon^{i} x_{i}$$
Any such dilatation is a group morphism and a Lie algebra morphism.

In fact the class of Carnot groups is characterised by the existence of dilatations.

\begin{prop}
Suppose that the Lie algebra $\mathfrak{g}$ admits an one parameter group
$\varepsilon \in (0,+\infty) \mapsto \delta_{\varepsilon}$ of simultaneously
diagonalisable Lie algebra isomorphisms. Then $\mathfrak{g}$ is the algebra of
a Carnot group.
\end{prop}

We can always find Euclidean inner products on $N$ such that the
decomposition $N \ = \ \sum_{i=1}^{m} V_{i}$ is an orthogonal sum. Let us pick
such an inner product and denote by $\| \cdot \|$ the Euclidean norm associated to
it.

We shall endow the group $N$ with a structure of a sub-Riemannian
manifold now. For this take the distribution obtained from left
translates of the space $V_{1}$. The metric on that distribution is
obtained by left translation of the inner product restricted to
$V_{1}$.

 The Carnot-Carath\'eodory  distance is 
$$d(x,y) \ = \ \inf \left\{ \int_{a}^{b} \| c^{-1} \dot{c} \| \mbox{
d}t \ \mbox{ : } \ c(a) = x , \ c(b) = y , \ c^{-1} \dot{c} \in
V_{1}
\right\}$$
The distance is obviously left invariant.

We collect the important facts to be known about Carnot groups:

\begin{enumerate}
\item[(a)] If $V_{1}$ Lie-generates the whole Lie algebra of $N$
then any two points can be joined by a horizontal path.
\item[(b)] The metric  topology and uniformity of $N$ are the same as Euclidean
topology and uniformity respective.
\item[(c)] The ball $B(0,r)$ looks roughly like the box
$\left\{ x \ = \ \sum_{i=1}^{m} x_{i} \ \mbox{ : }
\| x_{i} \| \leq r^{i} \right\}$.
\item[(d)] the Hausdorff measure $\mathcal{H}^{Q}$ is group
invariant and the Hausdorff dimension of a ball is $Q$.
\item[(e)] there is a one-parameter group of dilatations, where a
dilatation is an isomorphism $\delta_{\varepsilon}$ of $N$ which
transforms the distance $d$ in $\varepsilon d$.
\end{enumerate}

In Euclidean spaces, given $f: R^{n} \rightarrow R^{m}$ and
a fixed point $x \in R^{n}$, one considers the difference function:
$$X \in B(0,1) \subset R^{n} \  \mapsto \ \frac{f(x+ tX) - f(x)}{t} \in R^{m}$$
The convergence of the difference function as $t \rightarrow 0$ in
the uniform convergence gives rise to the concept of
differentiability in it's classical sense. The same convergence,
but in measure, leads to approximate differentiability. 
This and
another topologies might be considered (see Vodop'yanov
\cite{vodopis}, \cite{vodopis1}).

In the frame of Carnot groups the difference function can be written using only dilatations and the group operation. Indeed, for any function between Carnot groups
$f: G \rightarrow P$,  for  any fixed point $x \in G$ and $\varepsilon >0$  the finite difference function is defined by the formula:
$$X \in B(1) \subset G \  \mapsto \ \delta_{\varepsilon}^{-1} \left(f(x)^{-1}f\left(
x \delta_{\varepsilon}X\right) \right) \in P$$
In the expression of the finite difference function enters $\delta_{\varepsilon}^{-1}$ and $\delta_{\varepsilon}$, which are dilatations in $P$, respectively $G$.

Pansu's differentiability is obtained from uniform convergence of the difference
function when $\varepsilon \rightarrow 0$.

The derivative of a function $f: G \rightarrow P$ is linear in the sense
explained further.  For simplicity we shall consider only the case $G=P$. In this way we don't have to use a heavy notation for the dilatations.

\begin{defi}
Let $N$ be a Carnot group. The function
$F:N \rightarrow N$ is linear if
\begin{enumerate}
\item[(a)] $F$ is a {\it group} morphism,
\item[(b)] for any $\varepsilon > 0$ $F \circ \delta_{\varepsilon} \
= \ \delta_{\varepsilon} \circ F$.
\end{enumerate}
We shall denote by $HL(N)$ the group of invertible linear maps  of
$N$, called the  linear group of $N$.
\label{dlin}
\end{defi}

The condition (b) means that $F$, seen as an algebra morphism,
preserves the grading of $N$.

The definition of Pansu differentiability follows:

\begin{defi}
Let $f: N \rightarrow N$ and $x \in N$. We say that $f$ is
(Pansu) differentiable in the point $x$ if there is a linear
function $Df(x): N \rightarrow N$ such that
$$\sup \left\{ d(F_{\varepsilon}(y), Df(x)y) \ \mbox{ : } \ y \in B(0,1)
\right\}$$
converges to $0$ when $\varepsilon \rightarrow 0$. The functions $F_{\varepsilon}$
are the finite difference functions, defined by
$$F_{t} (y) \ = \ \delta_{t}^{-1} \left( f(x)^{-1} f(x
\delta_{t}y)\right)$$
\end{defi}
    
For the differentiability notion in a sub-Riemannian manifold the reader can consult Margulis, 
Mostow  \cite{marmos1} \cite{marmos2}, Vodop'yanov , Greshnov \cite{vodopis2},  \cite{vodopis3} or Buliga \cite{buliga2}. 
We shall use further the fact that isometries of a sub-Riemannian manifold are derivable in the 
sense of Pansu and the derivative is linear in the sense of the definition \ref{dlin}.
 
\section{Gromov-Hausdorff distance}
  The references for this section  are Gromov \cite{gromov}, chapter 3, and Burago \& al. \cite{burago} section 7.4.  There are several definitions of distances between metric spaces. The very fertile idea of introducing such distances belongs to Gromov.

In order to introduce the Hausdorff distance between metric spaces, recall
the Hausdorff distance between subsets of a metric space.

\begin{defi}
For any set $A \subset X$ of a metric space and any $\varepsilon > 0$ set
the $\varepsilon$ neighbourhood of $A$ to be
$$A_{\varepsilon} \ = \ \cup_{x \in A} B(x,\varepsilon)$$
The Hausdorff distance between $A,B \subset X$ is defined as
$$d_{H}^{X}(A,B) \ = \ \inf \left\{ \varepsilon > 0 \mbox{ : } A \subset B_{\varepsilon} \ , \ B \subset A_{\varepsilon} \right\}$$
\end{defi}

By considering all isometric embeddings of two metric spaces $X$, $Y$ into
an arbitrary metric space $Z$ we obtain the Hausdorff distance between $X$, $Y$ (Gromov \cite{gromov} definition 3.4).

\begin{defi}
The Hausdorff distance $d_{H}(X,Y)$ between metric spaces $X$ $Y$ is the infimum of the numbers
$$d_{H}^{Z}(f(X),g(Y))$$
for all isometric embeddings $f: X \rightarrow Z$, $g: Y \rightarrow Z$ in a
metric space $Z$.
\end{defi}

If $X$, $Y$ are compact then $d_{H}(X,Y) < + \infty$. Indeed, let $Z$ be the disjoint union of $X,Y$ and $M \ = \ \max \left\{ diam(X) , diam(Y) \right\}$.
Define the distance on $Z$ to be
$$d^{Z}(x,y) \ = \ \left\{ \begin{array}{ll}
d^{X}(x,y) & x, y \in X \\
d^{Y}(x,y) & x, y \in Y \\
\frac{1}{2} M & \mbox{ otherwise}
\end{array} \right. $$
Then $d_{H}^{Z}(X,Y) < + \infty$.

The Hausdorff distance between isometric spaces equals $0$. The converse is also true (Gromov {\it op. cit.} proposition 3.6) in the class of compact metric spaces.

\begin{thm}
If $X,Y$ are compact metric spaces such that $d_{H}(X,Y) = 0$ then
$X,Y$ are isometric.
\label{tgro}
\end{thm}

  \section{Deformations of a sub-Riemannian manifold}
  
 There are several deformations of a sub-Riemannian manifold around a point, which can 
 be studied in the in the metric spaces $CMS$ of isometry classes of compact metric spaces, 
 with the Gromov-Hausdorff  distance. For the isometry class of the metric space $(X,d)$ we shall 
 use he notation $[X,d]$. 
 
 The Ball-Box theorem (or the theorem concerning the existence of privileged frames) ensures us that 
 small closed balls in sub-Riemannian manifolds are compact. 
 
 \begin{defi}
 The metric profile associated to the 
 space $(M,d)$ is  the assignment (for small enough $\varepsilon > 0$) 
$$(x \in M, \varepsilon > 0) \ \mapsto \  \mathbb{P}^{m}_{\varepsilon}(x) = [\bar{B}(x,1) , \frac{1}{\varepsilon} d] \in CMS$$
\end{defi}

The celebrated Mitchell \cite{mit}  theorem 1 can be formulated in the following way: 

\begin{thm} (Mitchell, theorem1) The metric profile of a regular sub-Riemannian manifold can be prolonged by continuity in $\varepsilon = 0$. Moreover 
$$\mathbb{P}^{m}_{0}(p) = [\bar{B}(0,1), d_{N}]$$
the isometry class of the nilpotentization of the distribution at $p$.
\end{thm}

There are several proofs of this theorem. In order to understand them I shall introduce the notion of metric profile further. 

Consider   a privileged chart around $p \in M$. With this chart come the associated privileged 
frame, dilatations $\delta_{\varepsilon}$ and $\Delta_{\varepsilon}$. For any Riemannian 
metric $g$ one can define  deformations of this metric by the formula: 
$$g_{\varepsilon} (\Delta_{\varepsilon}X,\Delta_{\varepsilon}Y) = g(X,Y)$$

Let us begin by describing deformations induced from a privileged chart. These deformations are seen 
as curves in CMS, the space of isometry classes of compact metric spaces, with the Gromov-Hausdorff 
distance. 

The dilatation flow $\delta_{\varepsilon}$ induces the following deformation: 
$$[\delta, \varepsilon](t) = [\bar{B}(p,1), (\delta, \varepsilon)]$$
where  the distance $(\delta,\varepsilon)$ is given by 
$$(\delta,\varepsilon)(x,y) = d(\delta_{\varepsilon}^{-1}x, \delta_{\varepsilon}^{-1}y)$$
We can induce another deformation: let $(D_{\varepsilon}, g_{\varepsilon})$ be the pair  
distribution - metric on the distribution obtained by transport with $\delta_{\varepsilon}$, namely: 
$$D_{\varepsilon}(\delta_{\varepsilon}x) = D\delta_{\varepsilon}(x) D(x)$$
$$g_{\varepsilon}(\delta_{\varepsilon}x)(D\delta_{\varepsilon}(x)u,D\delta_{\varepsilon}(x)v) = 
g(x)(u,v)$$
for any $u,v \in T_{x}M$.  The deformation associated is 
$$[D, g,\varepsilon] = [ \bar{B}(p,1), (D, g,\varepsilon)]$$
where the notation $(D, g,\varepsilon)$ is used  for the CC distance in the 
sub-Riemannian manifold $(M,D_{\varepsilon}, g_{\varepsilon})$. 

It is an important remark that generally 
$$[\delta,\varepsilon] \not = [D, g,\varepsilon]$$
This is because the right-handed term is a distance given as the infimum of lengths of some horizontal curves. 

A sufficient condition for the equality to happen is that $\delta_{\varepsilon}$ has local convex data 
in the sense of the generalized Local-Global Principe to be found in Buliga \cite{buliga1} Section1.2.  In this case (small) CC balls in the manifold $(M,D_{\varepsilon}, g_{\varepsilon})$ are convex 
(that is there are geodesics connecting any two points in the closure of the ball, inside the closure of tha ball). This is not happening even if $M$ is a Carnot group. For example the closed balls in the Heisenberg group are not convex (with respect to the CC distance). 

Finally, a deformation is associated to the dilatations $\Delta_{\varepsilon}$ and pairs privileged frame - 
Riemannian metric $g$. This is simply

$$[\Delta,\varepsilon] = [\bar{B}(p,1), (\Delta,\varepsilon)]$$
where $(\Delta,\varepsilon)$ is the Riemannian distance induced by the Riemannian metric $g$ and 
the privileged frame. 

We shall see that all these deformations are particular  metric profiles. A general definition of a metric profile will be given further, after the discution of various proofs of Mitchell theorem 1.

\section{Mitchell theorem 1} 
 
 One can identify in the literature several proofs of this theorem. Exactly what is proven in each available variant of proof? The answer is: each proof basically shows that various metric profiles introduced in the previous section can be prolonged  to $0$. Each of this metric profiles are close in the GH distance 
 to the original metric profile of the CC distance. More precisely: 
 
 \begin{lema}
 Let $\mathbb{P}'(t)$ be any of the previously introduced metric profiles $[\delta, \varepsilon]$, 
 $[D,g,\varepsilon]$, $[\Delta, \varepsilon]$. Then 
 $$d_{GH}(\mathbb{P}^{m}_{\varepsilon}, \mathbb{P}'(\varepsilon)) = O(\varepsilon)$$
 \label{telema}
 \end{lema}
 
 The proof of this lemma reduces to a control problem. In the case of the profile $[\delta,\varepsilon]$, 
 this is Mitchell \cite{mit} lemma 1.2. 
      
  Mitchell \cite{mit} and Bella\"{\i}che \cite{bell} theorem 5.21, proposition 5.22,  proved the following: 
\begin{thm}
The metric profile $\varepsilon \mapsto [D,g,\varepsilon]$ can be prolonged at 
$\varepsilon = 0$ by continuity. We have  
$$[D,g,0]  = [\bar{B}(0,1), d_{N}]$$
\end{thm}

The corresponding  result of Gromov \cite{gromo} section 1.4B and Vodop'yanov \cite{vodopis3}  is: 
\begin{thm}
The metric profile $\varepsilon \mapsto [\Delta,\varepsilon]$ can be prolonged at 
$\varepsilon = 0$ by continuity. We have  
$$[\Delta,0]  = [\bar{B}(0,1), d_{N}]$$
\end{thm}

By the use of the Ball-Box theorem and any of the results before, one can obtain  an analogous result: 

\begin{thm}
The metric profile $\varepsilon \mapsto [\delta,\varepsilon]$ can be prolonged at 
$\varepsilon = 0$ by continuity. We have  
$$[\delta,0]  = [\bar{B}(0,1), d_{N}]$$
\end{thm}

The proofs of these theorems can be described as a manipulation of brackets associated with 
growth estimates in $\varepsilon$. 

Any of these theorems imply the Mitchell theorem 1, with the use of the approximation lemma \ref{telema}. But in fact these theorems are different statements in terms of metric profiles. 

I am interested to know if more information can be obtain from the use of a particular metric 
profile. We shall see that this is indeed the case. For example the curvature is a notion 
which is associated with a choice of such a profile. In the Riemannian case it does make no 
difference the choice of a metric profile. The phenomenon of dependence curvature --- metric profile 
is purely non-Riemannian. This path will not be pursued in this paper, where we shall use only the 
curvature given by the metric profile $\mathbb{P}^{m}$ associated with a metric space. 
        
\section{Metric profile}

The purpose of this section is two-folded. It serves as an introduction to the notion of metric profile. 
It is also written for further reference. For example the notion of approximate metric profile, useful in 
the understanding of the construction of a tangent bundle to a sub-Riemannian group (see Buliga 
\cite{buliga2}), will not be used in this paper. It will rather serve as an appetizer for the interested 
reader.

 The class of $\varepsilon$ nets (with arbitrary $\varepsilon$) 
in compact metric spaces will be denoted by NETS. In this paper nets always have positive 
separation. 

Likewise one can consider the classes 
$CMS_{a}$, $NETS_{a}$, of compact metric spaces (nets in compact metric spaces respectively) 
of diameter not greater than $a >0$.  The class $[NETS_{a}]$  with 
the Lipschitz distance is continuously embedded in $[CMS_{b}]$  with Gromov-Hausdorff distance, for 
any $b>a$.

We can define a notion of metric profile regardless to any distance. 

\begin{defi}
A metric profile is a curve $\mathbb{P}:[0,a] \rightarrow [CMS]$ such 
that:
\begin{enumerate}
\item[(a)] it is continuous at $0$,
\item[(b)]for any $b \in [0,a]$ and $\varepsilon \in (0,1]$ we have 
$$d_{GH} (\mathbb{P}(\varepsilon b), \mathbb{P}_{d_{b}}(\varepsilon,x)) - 
d_{GH}(\mathbb{P}(0),  \mathbb{P}_{d_{b}}(0,x)) \ 
\leq \ O(b) O(\varepsilon)$$
\end{enumerate}
We used here the notation $\mathbb{P}(b) = [\bar{B}(x,1),d_{b}]$ and $\mathbb{P}_{d_{b}}(\varepsilon,x) = [\bar{B}(x,1),\frac{1}{\varepsilon}d_{b}]$.
\label{dprofile}
\end{defi}

Note that in this definition is not stated that $\displaystyle \mathbb{P}(0) = \mathbb{P}_{d_{b}}(0)$. Look for example to the metric profile used by Gromov, namely $[\Delta, \varepsilon]$. For this profile we 
never have the mentioned equality, because $[\Delta,b](0)$ is always the Euclidean unit $n$ dimensional ball. Nevertheless this is a profile in the sense of the previous definition. 

\begin{defi}
The metric profile is nice if for all small enough $b$ we have $\displaystyle \mathbb{P}(0) = \mathbb{P}_{d_{b}}(0)$.
\end{defi}

The metric profile of a homogeneous space is just  a curve in the space 
$[CMS]$, continuous at $0$. Likewise, if we look at a homogeneous sub-Riemannian manifold, all 
metric profiles previously introduced are not depending on points in the manifold.

In order to give the definition of an approximate metric profile, we need a slightly modified 
version of  proposition 3.5, chapter 3, Gromov  \cite{gromov}. 

\begin{prop}
Let $(X_{i})_{i}$, $(Y_{i})_{i}$ be two sequences in $CMS$ such that 
$$d_{GH}(X_{i}, Y_{i}) \rightarrow 0$$
as $i \rightarrow \infty$. Then for any $\eta > 0$ and for any sequence $(N_{i})_{i} \subset NETS$
of $\eta$ nets $N_{i} \subset X_{i}$, there is a sequence $(M_{i})_{i} \subset NETS$ of 
$\eta + 2 d_{GH}(X_{i}, Y_{i}) + d_{GH}^{2}(X_{i}, Y_{i})$ nets $M_{i} \subset Y_{i}$ such that 
$$d_{Lip}(N_{i}, M_{i})  \ \leq \   2 d_{GH}(X_{i}, Y_{i}) + d_{GH}^{2}(X_{i}, Y_{i})$$
\end{prop}

\begin{cor}
Let $\mathbb{P}$ be a nice  metric profile, $\eta > 0$  and $\tilde{\mathbb{P}}_{\eta} : [0,1] \rightarrow [NETS]$ be  a curve such that $\tilde{\mathbb{P}}_{\eta}(a)$ is a $\eta$ net in $\mathbb{P}(a)$ for all $a$. 

Then there exists a function $\tilde{\mathbb{P}}_{\eta} : [0,1] \times [0,1] \rightarrow [NETS]$ such that 
\begin{enumerate}
\item[i)] $\tilde{\mathbb{P}}_{\eta}(a,1) = \tilde{\mathbb{P}}_{\eta}(a)$ for any $a$, 
\item[ii)] $\tilde{\mathbb{P}}_{\eta}(a,\varepsilon)$ is a $\eta + O(a)o(\varepsilon)$ net in 
$\mathbb{P}(a,\varepsilon)$, for all $a$, 
\item[iii)] the following estimate holds 
$$d_{Lip}(\tilde{\mathbb{P}}_{\eta}(a,1), \tilde{\mathbb{P}}_{\eta}(a,\varepsilon)) \ = \ 2\eta + 
O(a)O(\varepsilon)$$
\end{enumerate}
\end{cor}

The definition of an approximate metric profile follows. 

\begin{defi}
Let $\mathbb{P}$ be a nice metric profile. An approximate metric profile of $\mathbb{P}$ is a function 
$\tilde{\mathbb{P}}$ which satisfies the conclusions of the previous corollary, with the slight modification consisting in replacement of $\eta$ in the estimates by $O(\eta)$. 
\end{defi}

It would be  interesting to see what is happening in the case where an approximate metric profile is made by balls in discrete groups.

We shall define further the notion of curvature associated with a given metric profile. 

\begin{defi}
Suppose $\mathbb{P}$ is a nice  metric profile. Suppose moreover that 
$$d(\mathbb{P}(\varepsilon)(0), \mathbb{P}(\varepsilon)(a)) = O(\varepsilon a)$$
Then we shall call such a profile rectifiable at $\varepsilon = 0$. 

Two metric profiles $\mathbb{P}_{1}$, $\mathbb{P}_{2}$ which are rectifiable at $\varepsilon=0$ are 
equivalent if 
$$d(\mathbb(P_{1}(\varepsilon)(a), \mathbb{P}_{2}(\varepsilon)(a)) = o(a)$$
(for fixed $\varepsilon$). 

The curvature class of a metric profile $\mathbb{P}$ is the equivalence class of $\mathbb{P}$. 
\end{defi}

In particular cases we would like to be able to compute the curvature. This can be done by 
using homogeneous spaces.

\section{The homogeneous case}

To a homogeneous space we can associate the groups $Isom(X,d)$ and 
$Isom_{p}(X,d)$, of isometries (isometries wich fix  the point $p$ respectively) of 
$(X,d)$.  The coset class $Isom(X,d)/Isom_{p}(X,d)$ 
is homeomorphic with $(X,d)$ by the  map 
$$\pi: Isom(X,d)/Isom_{p}(X,d) \rightarrow X$$
The construction of the map $\pi$ is explained further. Let $p \in X$  be a fixed point. 
Pick a coset $f Isom_{p}(X,d)$ and define $\pi(f Isom_{p}(X,d)) = f(p)$. Obviously the definition is good. 

The inner action of $Isom_{p}(X,d)$ on $Isom(X,d)$ gives an action of $Isom_{p}(X,d)$ on the coset space $Isom(X,d)/Isom_{p}(X,d)$. This innes action is compatible with the action of $Isom_{p}(X,d)$ on 
$X$ in the sense: for any $h \in Isom_{p}(X,d)$ and for all $f \in Isom(X,d)$ we have 
$$ \pi(h f h^{-1} Isom_{p}(X,D)) = h (\pi( f Isom_{p}(X,d)))$$

In the case of a regular sub-Riemannian $(X,D,g)$ manifold we can associate to it the triple 
$(Isom(X,d), Isom_{p}(X,d), D, \gamma)$. 

The situation is as follows: note $G = Isom(X,d)$, $G_{0} = Isom_{p}(X,d)$. Then any right-invariant 
vectorfield on $G$ descends on a vectorfield on left cosets  $G/G_{0}$. In particular, if we endow 
$G$ with a right-invariant distribution, then $G/G_{0}$ is endowed with a distribution induced by the 
descent of any right invariant "horizontal" frame.  $G/G_{0}$ is not usually  a regular sub-Riemannian manifold. Look for example to the case: $G = H(1)$, $G_{0}$ is the one parameter group generated 
by an element of the distribution. Then $G/G_{0}$ is the Grushin plane, which is not a regular sub-Riemannian manifold. 

Consider on $G$ the right invariant distribution 
$$D" = Lie \  G_{0} + D'$$
such that $Lie \ G_{0} \cap D' = 0$ and $D'$ descends on the distribution $D$ on $G/G_{0}$ (if 
$G = Isom(X,d)$, $ G_{0} = Isom_{p}(X,d)$) . We have seen that the action of $G_{0}$ on 
$G/G_{0}$ which mimicks  the action of $Isom_{p}(X,d)$ on $(X,d)$ is the (descent of) the inner action. 
We are not wrong if we suppose that isometries preserve the distribution at $0$, which translates into the 
condition: for any $h \in G_{0}$ 
$$Ad_{h} D' \subset D'$$

We know one more thing about the homogeneous metric space $(X,d)$: its tangent cone.  Consider 
on $G$ with given distribution $D"$ the dilatations $\delta_{\varepsilon}$ and a privileged right-invariant 
basis around the neutral element. The knowledge of the tangent cone implies the following: 
\begin{enumerate}
\item[(a)]  we know some relations in the algebra $Lie \  G$, 
\item[(b)] we know that for any $h \in G_{0}$ $Ad_{h} \in HL(G,D")$, that is $Ad_{h}$ commutes with 
dilatations $\delta_{\varepsilon}$. 
\end{enumerate}

In conclusion, we can describe homogeneous metric spaces coming from sub-Riemannian manifolds
by looking to triples $(G,G_{0}, D")$, which satisfy certain relations. 

It goes without saying that we have also an Euclidean metric on the distribution $D"$. 

\begin{defi}

Let $(X,d)$ be a metric space and $p \in X$ a point such that the metric profile associated to $X,d)$ 
and $p$ can be prolonged at $\varepsilon = 0$ and it is rectifiable at $\varepsilon = 0$. We shall say 
that the curvature of $(X,d)$ at $p$ is $(G, G_{0}, D")$ if the metric profile at $p$ is equivalent with 
the metric profile of $G/G_{0}$ with respect to the distribution $D$ (the descent of $D"$). 

\end{defi}
 
 This definition is insinuating that $(G,G_{0}, D")$ (and the overlooked metric on $D"$) are uniquely 
 defined up to trivial transformation. We shall explore this issue in the final section, for a particular case. 
 
 As an exercise we want to compute all Riemannian homogeneous $n$ manifolds. So we are 
looking at groups $G$ which contain a subgroup $G_{0}$ such that: 
$$Lie \  G  \  = \ Lie \ G_{0} + D'$$ 
$$[Lie \ G_{0}, D'] \subset D'$$
$$[D',D'] \subset Lie \ G_{0}$$
and for any $x \in Lie \  G_{0}$ the restriction of $ad_{x}$ on $D'$ is antisymmetric . Moreover, $D'$ has dimension $n$. 
For example, when $n=2$ we have two cases. The first case is $G$ 3 dimensional, with a basis 
$X_{0}, X_{1}, X_{2}$ for $Lie \ G$, such that $X_{0}$ generates $Lie \ G_{0}$ and $X_{1}, X_{2}$ 
generate $D'$. The bracket relations that we know are: 
$$[X_{0}, X_{1}] = a X_{1} + b X_{2}$$
$$[X_{0}, X_{2}] = c X_{1} + d X_{2}$$
$$[X_{1}, X_{2}] = e X_{0}$$
From Jacobi identity we get $e(a+d) = 0$ and from the condition $ad_{X_{0}}$ restricted to $D'$ antisymmetric we get $a = 0$, $d = 0$, $b+c = 0 $. Therefore we have
$$[X_{0}, X_{1}] =  b X_{2}$$
$$[X_{0}, X_{2}] = -b X_{1}$$
$$[X_{1}, X_{2}] = e X_{0}$$
 We have a one dimensional family of homogeneous Riemannian surfaces, where the curvature can 
 be measured by $b/e$(except $e=0$, but the factor space  is the Euclidean plane; see also next case). 
 
 The second case is $dim \ G_{0} = 0$ and $G$ is abelian 2 dimensional. But this is trivial, moreover, it is 
 contained in the previous case. 
 
 This is well known and seems to be related to the Cartan method of equivalence.

\section{Application: curvature of contact 3 manifolds}

  Contact manifolds are particular cases of sub-Riemannian manifolds. The contact distribution 
  is completely non-integrable. By using natural normalization of the contact form (see for example 
  Bieliavski, Flbel, Gorodski \cite{biefago} or Hughen \cite{hughen}) we can uniquely associate to a contact structure, endowed with a metric on the contact distribution, a sub-Riemannian manifold. The nilpotentization  of the contact  distribution is 
  always a Heisenberg group. 
  
  The horizontal linear maps os the Heisenberg group are known. Moreover, the group of isometries 
  of $H(n))$ which preserve the origin is   $SU(n)$. 
  
  If we want to look for all homogeneous contact 3 manifolds, we have to consider two cases. 
  The first case is $G$ 4 dimensional, with a basis for $Lie \  G $ given by $X_{0}, X_{1}, X_{2}, X_{3}$, such that $X_{0}$ is a basis for $Lie \ G_{0}$, 
$$[X_{1}, X_{2}] = X_{3}+ a X_{0}$$
$$[X_{1}, X_{3}] = b X_{0} + A X_{1} + B X_{2}$$
$$[X_{2}, X_{3}] = c X_{0} + C X_{1} + D X_{2}$$ 
(which comes from the knowledge of the nilpotentization and from the condition $[D',D'] \subset Lie \ G_{0}$), 
$$[X_{0}, X_{1}] = d X_{1} + e X_{2}$$
$$[X_{0}, X_{2}] = f X_{1} + g X_{2}$$
and the $ad_{X_{0}}$ condition is $d=g =0$,  $e + f = 0$ and $[X_{0}, X_{3}] = 0$.

    We have to use further the Jacobi identities. We begin with: 
$$[X_{1}, [X_{2}, X_{3}]] + [X_{2}, [X_{3}, X_{1}]] +[X_{3}, [X_{1}, X_{2}]] = 0$$
This gives:  $A+D = 0$, $eC = 0$, $e b = 0$. 

The next relation is: 
$$[X_{0}, [X_{2}, X_{3}]] + [X_{2}, [X_{3}, X_{0}]] +[X_{3}, [X_{0}, X_{2}]] = 0$$
This gives: $e(A-D) = 0$, $e(B+C)= 0$. 

The relation: 
$$[X_{0}, [X_{1}, X_{2}]] + [X_{1}, [X_{2}, X_{0}]] +[X_{2}, [X_{0}, X_{1}]] = 0$$
gives nothing new. 

We continue with: 
$$[X_{0}, [X_{1}, X_{3}]] + [X_{1}, [X_{3}, X_{0}]] +[X_{3}, [X_{0}, X_{1}]] = 0$$
which lead to nothing new. 

If $e=0$ then we have $f=0$,  $A+D=0$, and we get the relations: 
$$[X_{1}, X_{2}] = X_{3}+ a X_{0}$$
$$[X_{1}, X_{3}] = b X_{0} + A X_{1} + B X_{2}$$
$$[X_{2}, X_{3}] = c X_{0} + C X_{1} -A X_{2}$$
$$[X_{0}, X_{1}] = 0$$
$$[X_{0}, X_{2}] = 0$$
$$[X_{0}, X_{3}] = 0$$
By a change of basis: $X_{3}' = X_{3} + a X_{0}$, ... , we arrive to the description of $G$ as a direct 
sum of a 3 dimensional group with $G_{0} = S(1)$. This is in reality a singular case (in the sense that 
$G_{0}$ is not really needed in the construction: it is added and after factorized out without any consequences).

If $e \not = 0$ then we have the relations
$$[X_{1}, X_{2}] = X_{3}+ a X_{0}$$
$$[X_{1}, X_{3}] = 0$$
$$[X_{2}, X_{3}] = c X_{0} $$
$$[X_{0}, X_{1}] = e X_{2}$$
$$[X_{0}, X_{2}] = -eX_{1}$$
$$[X_{0}, X_{3}] = 0$$

These form a 2 dimensional family of  homogeneous spaces which are not groups. 

The second case is $dim \ G_{0}  = 0$ and $G$ is 3  dimensional, with a basis for $Lie \  G $ given by 
$ X_{1}, X_{2}, X_{3}$, such that: 
$$[X_{1}, X_{2}] = X_{3}$$
$$[X_{1}, X_{3}] =  A X_{1} + B X_{2}$$
$$[X_{2}, X_{3}] =  C X_{1} + D X_{2}$$ 
The Jacobi identity 
$$[X_{1}, [X_{2}, X_{3}]] + [X_{2}, [X_{3}, X_{1}]] +[X_{3}, [X_{1}, X_{2}]] = 0$$
gives the relation $A+D = 0$, therefore we recover  previous case. 

We have a 3 dimensional family of regular homogeneous spaces which are also groups. Particular examples are: 
SO(3), SL(2,R), E(1,1). 

\section{Classification of curvatures}

We shall prove in this section that for any two 3 dimensional homogeneous spaces which are also 
groups, if they have the same curvature class then they are isometric. This will partially solve the problem of classification of curvatures for 3 dimensional contact manifolds. 

More specifically we shall prove the following: 

\begin{thm}
Let $G_{1}$, $G_{2}$ be two 3 dimensional groups. 
We identify the Lie algebras and hence we have two brackets on $\mathbb{R}^{3}$ denoted by 
$[ \cdot , \cdot ]_{i}$, $i = 1,2$. 

 Define also $\delta_{\varepsilon}(X_{1}) = \varepsilon X_{1}$, 
 $\delta_{\varepsilon}(X_{2}) = \varepsilon X_{2}$ and 
 $\delta_{\varepsilon}(X_{3}) = \varepsilon^{2} X_{3}$.

Let $d_{1}$, $d_{2}$ be the CC distances on $G_{1}$, $G_{2}$ with respect to the left invariant 
distributions generated by $X_{1}, X_{2}$, transported on (a neighbourhood of $0$ of) 
$\mathbb{R}^{3}$. 

 Suppose that we have  the bracket relations: 
$$[X_{1}, X_{2}]_{i} = X_{3}$$
$$[X_{1}, X_{3}]_{i} =  A_{i} X_{1} + B_{i} X_{2}$$
$$[X_{2}, X_{3}]_{i} =  -B_{i} X_{1} + D_{i} X_{2}$$ 

If $d_{1}(\delta_{\varepsilon} x, \delta_{\varepsilon} y) - 
d_{1}(\delta_{\varepsilon} x, \delta_{\varepsilon} y) = 0(\varepsilon^{2})$ uniformly with respect 
to $x, y$ in a compact neighbourhood of $0$, then the Lie brackets are identical. 
\end{thm}
 
The proof uses the Baker-Campbell-Hausdorff formula and the Ball Box theorem. The hypothesis implies that 
$$d_{1}^{2}(\delta_{\varepsilon} x, \delta_{\varepsilon} y) - 
d_{1}^{2}(\delta_{\varepsilon} x, \delta_{\varepsilon} y) = 0(\varepsilon^{4})$$
Each distance $d_{i}$ is left invariant. We shall note 
$$\| u \|_{i} = d_{i}(0,u)$$
We know from the Ball Box theorem that $\| u \|^{2}$ is comparable with $\| u_{1} \|^{2} + \mid u_{2} \mid$, where $u= u_{1} + u_{2}$ is the decomposition of $u$ into the horizontal part 
$u_{1} \in \ span \ \left\{ X_{1}, X_{2} \right\}$ and the vertical part $u_{2} \in \ span \ \left\{ X_{3} \right\}$. 

We shall denote by $\opunu$, $\opdoi$ the operations in $G_{1}$, $G_{2}$ respectively.  The hypothesis becomes: 
$$\ | \delta_{\varepsilon} (-x) \opunu \delta_{\varepsilon}y \|^{2}_{1} - 
\ | \delta_{\varepsilon} (-x) \opdoi \delta_{\varepsilon}y \|^{2}_{2} = 0(\varepsilon^{4})$$
From the Baker-Campbell-Hausdorff formula and the bracket relations we see that we can approximate 
  $ \delta_{\varepsilon} (-x) \opunu \delta_{\varepsilon}y$ up to $o(\varepsilon^{4})$ by using only terms 
  in the Baker-Campbell-Hausdorff formula which contain at most two brackets. Same is true for the operation $\opdoi$. 
  
  Moreover the norms $\| \cdot \|_{i}$ can be estimated from the Ball Box theorem. From here a careful  
  computation resumes the proof.

\end{document}